\documentclass[12pt,a4paper]{article}

\usepackage[utf8]{inputenc}
\usepackage[T1]{fontenc}
\usepackage{lmodern}
\usepackage{amsmath,amsfonts,amssymb,amsthm}
\usepackage{bbm}
\newcommand*{\indic}{\ensuremath{\mathbbm{1}}}
\usepackage{graphicx}
\usepackage{tikz} 
\DeclareMathOperator*{\argmin}{argmin}
\DeclareMathOperator*{\arcccsh}{ArcCsch}
\DeclareMathOperator*{\arcsinh}{ArcSinh}
\title{One step futher: an explicit solution to Robbins' problem when $n=4$} \author{Rémi Dendievel\footnote{Mathematics
    Department, Universit\'e libre de Bruxelles, Belgium}, \quad Yvik
  Swan\footnote{Mathematics Department, Universit\'e de Li\`ege,
    Belgium}} \date{}

\newcommand{\domaine}[6]{\begin{tikzpicture}[scale=#6]
  \draw (0,0) rectangle (1,1);
  \coordinate (q1) at (1/4,1/4);
  \coordinate (q2) at (1/2,1/2);
  \coordinate (q3) at (3/4,3/4);
  \filldraw[fill=black!3] (0,1/3) -- (q1) -- (1/3,0) -- (0,0) --
  cycle;
  \filldraw[fill=black!10] (0,1/3) -- (q1) -- (1/3,0) -- (2/3,0) --
  (q2) -- (0,2/3) -- cycle;
  \filldraw[fill=black!15] (0,2/3) -- (q2) -- (1/3,1) -- (0,1) --
  cycle;
  \filldraw[fill=black!20] (2/3,0) -- (q2) -- (1,1/3) -- (1,0) --
  cycle;
  \filldraw[fill=black!25] (1/3,1) -- (2/3,1) -- (q3) -- (1,2/3) --
  (1,1/3) -- (q2) -- cycle;
  \filldraw[fill=black!30] (2/3,1)--(1,1)--(1,2/3)--(q3)--cycle;

  \node at (0.12,0.12) {#1};
  \node at (0.35,0.35) {#2};
  \node at (0.25,0.8) {#3};
  \node at (0.8,0.25) {#3};
  \node at (0.65,0.65) {#4};
  \node at (0.88,0.88) {#5};
  \node[below] at (1/2, -0.2) {$x_1$};
  \node[left] at (-0.2, 1/2) {$x_2$};
  \node[below, left] at (-0.02,-0.02) {\scriptsize$0$};
  \node[below] at (1,0) {\scriptsize$1$};
  \node[left] at (0,1) {\scriptsize$1$};
  \node[left] at (0,1/3) {\scriptsize$1/3$};
  \node[left] at (0,2/3) {\scriptsize$2/3$};
  \node[below] at (1/3,0){\scriptsize$1/3$};
  \node[below] at (2/3,0){\scriptsize$2/3$};

  \draw[] (0,1/3) -- (q1) -- (1/3,0) -- (0,0) --
  cycle;
  \draw[] (0,1/3) -- (q1) -- (1/3,0) -- (2/3,0) --
  (q2) -- (0,2/3) -- cycle;
  \draw[] (0,2/3) -- (q2) -- (1/3,1) -- (0,1) --
  cycle;
  \draw[] (2/3,0) -- (q2) -- (1,1/3) -- (1,0) --
  cycle;
  \draw[] (1/3,1) -- (2/3,1) -- (q3) -- (1,2/3) --
  (1,1/3) -- (q2) -- cycle;
  \draw[] (2/3,1)--(1,1)--(1,2/3)--(q3)--cycle;

\end{tikzpicture}}

\begin{document}
\maketitle{}

\begin{abstract}
  Fix some $n \in \mathbb{N}$ and let $X_1, X_2,\dots, X_n$ be
  independent random variables drawn from the uniform distribution on
  $[0,1]$. A decision maker is shown the variables sequentially and,
  after each observation, must decide whether or not to keep the
  current one, with payoff the overall rank of the selected
  observation.  Decisions are final: no recall is allowed, no regret
  is tolerated. The objective is to act in such a way as to minimize
  the expected payoff.  In this note we give the explicit solution to
  this problem, known as Robbins' problem of optimal stopping, when
  $n=4$.
\end{abstract}

\section{Introduction}
\label{sec:introduction-history}
Robbins' problem (of optimal stopping) consists in studying the
mathematical properties of the optimal strategy in the following
sequential selection problem.

\begin{quote}
  \emph{Fix some $n \in \mathbb{N}$ and let $X_1, X_2,\dots, X_n$ be
    independent random variables drawn from the uniform distribution
    on $[0,1]$. A decision maker is shown the variables sequentially
    and, after each observation, must decide whether or not to keep
    the current one. The payoff is $R_k$, the overall rank of the selected
    observation, with the convention 
    \begin{equation*}
      R_k = \sum_{i=1}^n \mathbb{I}(X_i \le X_k)
    \end{equation*}
    (and $\mathbb{I}(A)$ the indicator function of $A$).  Decisions
    are final: no recall is allowed, no regret is tolerated. The total
    number of observations is known to the decision maker.  The
    objective is to act in such a way as to minimize the expected
    overal rank of the selected observation. }
\end{quote}

In the sequel we use the shorthand $RP(n)$ to refer to the above
problem with $n$ arrivals. Solving Robbins' problem consists in
describing $\tau_n^{\star}$, the optimal stopping rule, computing
$v(n)$, the optimal expected rank obtainable with $n$ observations,
understanding the main traits of $\tau^{\star}_n$ as $n$ grows large
and obtaining the limiting value $\lim_{n\to \infty} v(n)=v$. Coaxed
by Prof.\ Herbert Robbins in the early 1990's (see Bruss 2005),
several independent teams devoted a significant amount of effort on
this seemingly innocuous problem. All have come to the conclusion that
the problem is \emph{``very hard''}. So much so that a complete
solution to Robbins' problem still eludes us to this date.

Robbins and coauthors (see Chow \emph{et al.} 1964) solve a
\emph{no-information} version of the problem, in which the decision
maker is not given the \emph{values} of the observations but only
their \emph{relative ranks}. Denoting $W(n)$ the corresponding
expected rank, Chow \emph{et al.} (1964) provide the
optimal strategy and manage an analytic \emph{tour de force} to prove
that $ W(n) \to W \approx 3.8695$, as $n \to \infty$.  Clearly $W(n)
\ge v(n)$ for all $n \ge1$, and hence we deduce that
\begin{displaymath}
  v \le 3.8695.
\end{displaymath}
Of course the full-information $RP(n)$ is much more favorable to the
decision maker and we thus expect $v(n)$ and $v$ to be, in fact, much
smaller than $W(n)$ and $W$, respectively.

Taking advantage of the knowledge of the values of the arrivals it is
natural to consider the class of stopping rules of the form
\begin{equation}\label{eq:2}
  \tau^{(n)} = \inf\left\{ k \ge 1 \, | \, X_k \le c_k^{(n)} \right\},
\end{equation}  
which we will call \emph{memoryless threshold rules}.  Bruss and
Ferguson (1996) prove that there exists a unique optimal
sequence (that is, optimal among memoryless threshold rules) which is
stepwise increasing in $n$. Also it is shown in Assaf and Samuel-Cahn
(1996) and in Bruss and Ferguson (1993) that
if $\tau^{(n)}$ is given by a sequence of increasing thresholds $0 <
a_1 \le a_2 \le \ldots \le a_n=1$, then
  \begin{align*}
    E \left( R_{\tau^{(n)}} \right) & =1+\dfrac{1}{2}{{
        \displaystyle{\sum_{k=1}^{n-1}}(n-k)}a_k^2{
        \displaystyle{\prod_{j=1}^{k-1}}(1-a_j)}} +\dfrac{1}{2}{
      \displaystyle{\sum_{k=1}^{n}{\prod_{j=1}^{k-1}(1-a_j)}}}
    \displaystyle{{\sum_{j=1}^{k-1}}\frac{(a_k-a_j)^2}{1-a_j}}
    \end{align*}
    with $R_{\tau^{(n)}}$ the rank of the observation selected by
    applying the stopping rule $\tau^{(n)}$. Clearly $v(n) \le E
    \left( R_{\tau^{(n)}} \right)$ for all $n$.  It is straightforward
    to optimize this expression over all possible thresholds (at least
    numerically) to obtain the values for $V(n) = \inf_{\tau^{(n)}} E
    \left( R_{\tau^{(n)}} \right)$ reported in Table~\ref{tab:1}.
\begin{table}[!]
  \centering
    \begin{equation*}
      \begin{array}{c|ccccccc}
n & 1 & 2 & 3 & 4 & 5 & 20 & 50 \\
\hline
V & 1 & 1.25 & 1.4009 & 1.5065 & 1.5861 &  1.9890 & 2.1482
      \end{array}
    \end{equation*}
\caption{Values of the memoryless optimal expected rank}\label{tab:1}
\end{table}
See Bruss and Ferguson (1996, Table 1b) (up to a minor correction of a typo for
their $V(4)$) or Bruss and Ferguson (1993) where the
computations are pushed as far as the case $n=800$.  Assaf and
Samuel-Cahn (1996) further explore rules based on
suboptimal thresholds of the form $ a_{k}^{(n)} = {\sum_{j=0}^mc_j
  k^j}/({n-k+c}) \wedge 1$ and mention numerical computations showing
that for $m=2$ the optimal coefficients are $c_0=1.77$, $c_1 = 0.54$
and $c_2 = -0.27$ yielding $V = \lim_{n\to \infty}V(n) \le
2.3268\cdots$ (our conclusion is slightly different to their value
2.3267; this is perhaps due to rounding errors in their computation)
and therefore
\begin{equation*}
  v \le  2.3268
\end{equation*}
(which is already an important improvement on the optimal
no-information value).  Although we still do not know the exact value
of $V$, Bruss and Ferguson (1993) extrapolate $V =
2.32659$ and Assaf and Samuel-Cahn (1996) prove that $V
\ge2.29558$, hence not much improvement on $v$ can be hoped for by
further exploring memoryless threshold rules of the form \eqref{eq:2}.

Intriguingly we know that there must exist rules which provide strict
improvement on those of the form \eqref{eq:2} because Bruss and
Ferguson (1993) prove that $v(n) < V(n)$ for all $n\ge1$,
i.e.\ even the optimal memoryless rule is strictly sub-optimal at
every $n$ for $RP(n)$. Meier and S\"ogner (2014) study
variations on the memoryless threshold rules wherein relative ranks
are taken into account and manage to lower the upper bound to obtain
an expected rank of 2.31301. This improvement is, however, not
significant enough even to answer whether or not $v$ is strictly
smaller than $V$ or not.

Several authors (e.g.\ Gnedin 2007, Bruss and Swan 2009 and
Gnedin and Iksanov 2011) have considered an alternative approach to
Robbins's problem by embedding it in a Poisson process.  Gnedin (2007)
proves that the memoryless stopping rules remain sub-optimal even in a
Poisson limiting model, i.e.\ there must exist stopping rules which
take the history of the arrival process into account and which provide
a strict improvement (even in a Poissonian limit) on the optimal
memoryless threshold rule. As can be seen from Bruss and Swan (2009),
embedding the problem in a Poisson arrival process yields several
advantages and opens several new veins of research on this fascinating
problem (see also Gnedin and Iksanov 2011) but still does not provide
satisfactory solutions to the original problem.

Backward induction guarantees the existence of an optimal strategy
$\tau_{\star}^{(n)}$ and provides, in principle, a way to compute
it. Hence for each $n\ge1$ there must exist threshold functions
$ h_k^{(n)}: [0,1]^{k-1} \to [0,1], k=1, \ldots, n-1$ such that
the optimal stopping rule is
\begin{equation*}
\tau_{\star}^{(n)} = \inf \left\{ k \, | \, X_k \le  h_k^{(n)}(X_1,
   \ldots, X_{k-1}) \right\}. 
\end{equation*}
Bruss and Ferguson (1993, 1996) prove that the threshold functions are
pointwise increasing but depend in a {non-monotone way} on {all} the
values of the previous arrivals and any loss of information results in
the loss of optimality.  This last point is referred to as \emph{full
  history dependence} of the optimal policy.  A consequence is that
any direct computations related to this optimal strategy are
fiendishly complicated and even computer simulations with modern-day
technology cannot bring any intuition even for moderate values of $n$
(double exponential complexity). We refer the reader to Bruss (2005)
for further information on the problem and its history.

To this date the optimal policy was only explicitly known in the case
$n=2$ (basically trivial) and $n=3$ (provided by Assaf and Samuel-Cahn
1996), with values $v(2) = 1.25 \mbox{ and } v(3) = 1.3915\cdots$,
respectively. The purpose of this note is to provide a modest
complement to the literature by solving the case $n=4$.  We will
derive the optimal threshold functions $h_1^{(4)}$, $h_2^{(4)}(x_1)$
and $h_3^{(4)}(x_1, x_2)$, whose behaviour is a complicated function
of the past data, see Section \ref{sec:expl-solut-case} for details)
and compute the value $ v(4) = 1.4932\cdots$ which is remarkably close
to the optimal memoryless value $V(4) = 1.5065$ from
Table~\ref{tab:1}. For the sake of completeness we also provide a
proof for the optimal strategies and values in the cases $n=2$ and
$n=3$.  As far as we can see there is no easy way to generalize our
result to higher values of $n$.

\section{Solution for the cases $n=2$ and $n=3$}
\label{sec:expl-solut-cases}

The case $n=2$ is nearly trivial. Indeed the threshold value at step~2
must be taken as $1$, and only $h_1$ needs to be computed (here and
throughout we drop the superscript $(n)$ for the thresholds). Define
$G(h)$ as the expected {rank of the selected value by using a strategy
  with threshold $h_1=h$}.  This expression is minimal for $h_1=1/2$
and we immediately conclude $v(2)=5/4$ (which is obviously the same
value as $V(2)$ in Table~\ref{tab:1}).

We now tackle the case $n=3$. We know that $h_3 = 1$ and must
determine the thresholds $h_1$ and $h_2(x_1)$. Define, in the same
fashion as above, $G_{x_1}(h)$ as ``the expected rank of the selected
variable given $X_1=x_1$ if we start to play at step~2 by using a
threshold value set to $h$''. Direct computations yield
\begin{equation}
  \label{eq:compus-G-case-n-3}
  G_{x_1}(h) = \frac32 + h^2 - h + (1-x_1)(1-h) + (h-x_1)_+,
\end{equation}
where $y_+ = \max(y, 0)$.
\begin{figure}[h]
  \centering
  \begin{minipage}{3cm}\centering
    \textsc{Case $A_1$}\\[2pt]\hrule
  \begin{tikzpicture}[scale=2]
    \draw[>=stealth,->] (0,0) -- (1.1,0);
    \draw[>=stealth,->] (0,0) -- (0,1.1);
    \node[above] at (0,1.1) {\footnotesize$G_{x_1}(h)$};
    \node[below,right] at (1.1,0) {\footnotesize$h$};
    \coordinate (a) at (0.75,0.22);
    \draw[thin] (0.1,0.9) to[out=280,in=220] (a)
    to[out=65,in=265] (0.95,0.9);
    \draw[dotted] (a) -- (0.75,0) node[below] {\scriptsize$x_1$};
  \end{tikzpicture}
  \end{minipage}\qquad
  \begin{minipage}{3cm}\centering
    \textsc{Case $A_2$}\\[2pt]\hrule
    \begin{tikzpicture}[scale=2]
    \draw[>=stealth,->] (0,0) -- (1.1,0);
    \draw[>=stealth,->] (0,0) -- (0,1.1);
    \node[above] at (0,1.1) {\footnotesize$G_{x_1}(h)$};
    \node[below,right] at (1.1,0) {\footnotesize$h$};
    \coordinate (a) at (0.6,0.22);
    \draw[thin] (0.1,0.9) to[out=280,in=170] (a)
    to[out=25,in=265] (0.95,0.9);
    \draw[dotted] (a) -- (0.6,0) node[below] {\scriptsize$x_1$};
  \end{tikzpicture}
  \end{minipage}\qquad
  \begin{minipage}{3cm}\centering
    \textsc{Case $A_3$}\\[2pt]\hrule
    \begin{tikzpicture}[scale=2]
    \draw[>=stealth,->] (0,0) -- (1.1,0);
    \draw[>=stealth,->] (0,0) -- (0,1.1);
    \node[above] at (0,1.1) {\footnotesize$G_{x_1}(h)$};
    \node[below,right] at (1.1,0) {\footnotesize$h$};
    \coordinate (a) at (0.25,0.22);
    \draw[thin] (0.05,0.9) to[out=275,in=120] (a)
    to[out=325,in=265] (0.9,0.9);
    \draw[dotted] (a) -- (0.25,0) node[below] {\scriptsize$x_1$};
  \end{tikzpicture}
  \end{minipage}
  \caption{The three generic situations we must study in order to find
    the expression of the minimizer of $G_{x_1}$.}
  \label{fig:ABC}
\end{figure}
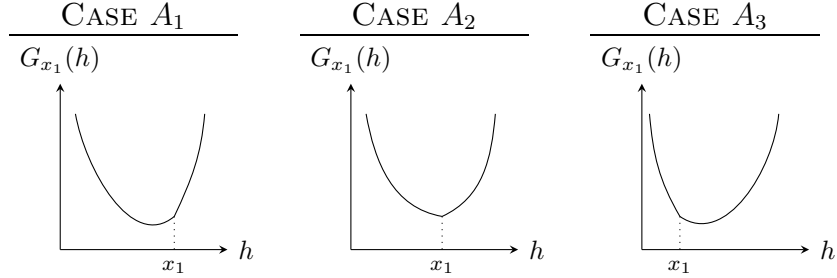

Minimizing in $h$ this expected rank we find that we must distinguish
three cases (see Figure~\ref{fig:ABC}) to get
\begin{equation}
  \label{eq:arminG-case-n-3}
  \argmin_{h\in[0,1]} G_{x_1}(h) = \begin{cases}
    \tfrac{1-x_1}{2} & \text{ if \:} 0 \le x_1 < \tfrac13 \text{ (case $A_1$)}\\
    x_1 & \text{ if \:} \tfrac13 \le x_1 < \tfrac23 \text{ (case $A_2$)}\\
    1-\tfrac{x_1}{2} & \text{ if \:} \tfrac23 \le x_1 \le 1 \text{ (case
      $A_3$)}
  \end{cases},
\end{equation}
from which we deduce $h_2(x_1)$, the optimal threshold at step 2.

By the optimality principle, the value of the threshold $h_1$ must be
a solution to the indifference equation
\begin{equation}
  \label{eq:indiff-h1-case-n-3}
  1 + 2h_1 = G_{h_1}(h_2(h_1))
\end{equation}
(i.e.\ the expected rank for choosing an arrival with value $h_1$ is
the same as for continuing and acting optimally thereafter).
Solutions of \eqref{eq:indiff-h1-case-n-3} are outside of $[0,1]$ both
when $h_1<1/3$ (case $A_1$) and $2/3\le h_1\le 1$ (case $A_3$). In
situation $A_2$ the equation becomes
\begin{equation*} 
  \label{eq:indiff-explicit-n-3}
  1 +2h_1 = \frac32 + h_1^2 - h + (1-h_1)^2,
\end{equation*}
with solution $h_1 = (5-\sqrt{13})/4$. This leads to 
the same conclusion as Assaf and Samuel-Cahn (1996),
namely that the optimal thresholds for $RP(3)$ are
\begin{equation*}
  h_1 = \frac{5-\sqrt{13}}{4}, \quad 
  h_2(x_1) =
  \begin{cases}
    x_1 & \text{ if \:}  h_1  \le x_1 \le \tfrac23\\
    1-x_1/2 & \text{ if \:} \tfrac23\le x_1 \le 1
  \end{cases}
\end{equation*}
(and $h_3 = 1$) 
providing us with the value  
\begin{equation*}
  v(3)  =  \frac{341}{144} -\frac{13}{48}  \sqrt{13} = 1.39155\cdots
\end{equation*}
which is remarkably close to the corresponding memoryless value $V(3)$ in
Table~\ref{tab:1}.

\section{Solution for the case $n=4$}
\label{sec:expl-solut-case}

As anticipated, in this section we prove the main contribution of this
note, namely 
\begin{equation}
  \label{eq:1}
  v(4) = 1.4932\cdots.
\end{equation}


The dynamic programming approach requires to find the optimal behaviour
at some specific step $k$ given a length $k-1$ history, by letting $k$
go backwards from $n$ to $1$. 
Our plan is thus simple : we start considering the best action at time
$k=4$, then we proceed backwards and end with the case $k=1$. For each
$k$, we fix a history $X_1=x_1$, $X_2=x_2$,\dots,
$X_{k-1}=x_{k-1}$. We know from Bruss and Ferguson (1993)
that the optimal action is defined by a threshold $h_k(x_1, \ldots,
x_{k-1})$: keep $X_k$ if less than $h_k(x_1, \ldots, x_{k-1})$,
otherwise discard it. Our purpose is to determine the exact
expressions for $h_k(x_1, \ldots, x_{k-1}), k=1, 2, 3, 4$.
\medskip

\noindent\textbf{Step 4.\:} Suppose that $(X_1, X_2, X_3)
=(x_1, x_2, x_3)$
has been observed and we only enter the game at step $4$ before
learning the value of $X_4$. Since this is the last step, we must
accept it whatever its value may be. This is the optimal behaviour,
and $h_4(x_1,x_2,x_3) = 1$, for all $(x_1,x_2,x_3)\in [0,1]^3$.
\medskip

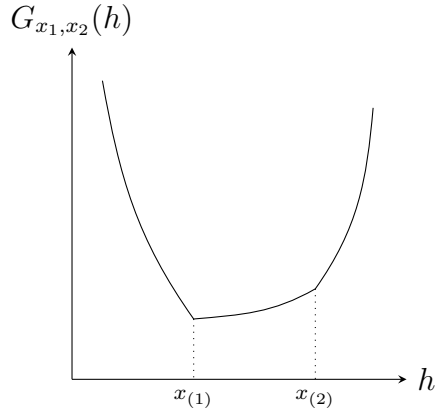
\begin{figure}[h]
  \centering
  \begin{tikzpicture}[scale=4]
    \draw[>=stealth,->] (0,0) -- (1.1,0);
    \draw[>=stealth,->] (0,0) -- (0,1.1);
    \node[above] at (0,1.1) {$G_{x_1,x_2}(h)$};
    \node[below,right] at (1.1,0) {$h$};
    \coordinate (a1) at (0.4,0.2);
    \coordinate (a2) at (0.8,0.3);
    \draw[thin] (0.1,0.99) to[out=280,in=125]
    (a1) to[out=5,in=210]
    (a2) to[out=55,in=265] (0.99,0.90);
    \draw[dotted] (a1) -- (0.4,0) node[below] {\scriptsize$x_{(1)}$};
    \draw[dotted] (a2) -- (0.8,0) node[below] {\scriptsize$x_{(2)}$};
  \end{tikzpicture}
  \caption{Graph of $G_{x_1,x_2}(\cdot)$ for one particular
    history. As in the case $n=3$, the minimum will be given by the
    minimizer of one of the parabolae or by one of the past
    observations. In our case ($n=4$), this leads to $5$ cases.}
  \label{fig:graphg}
\end{figure}

\begin{figure}[h]
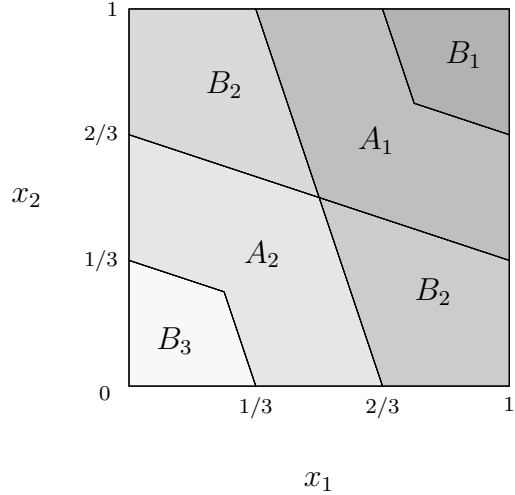

  \centering
  \domaine{$B_3$}{$A_2$}{$B_2$}{$A_1$}{$B_1$}{5}
  \caption{The regions $A_1,A_2,B_1,B_2,B_3$ are circumscribed by the
    borders of $[0,1]^2$ and the lines $x_2=(3-x_1)/3$,
    $x_2=(2-x_1)/3$, $x_2=(1-x_1)/3$, $x_2=3-3x_1$, $x_2=2-3x_1$,
    $x_2=1-3x_1$.}
  \label{fig:h3-regions}
\end{figure}

\noindent\textbf{Step 3.\:} Suppose that $(X_1, X_2)=(x_1, x_2)$ has
been observed and we enter the game at step $3$ before learning the
value of $X_3$. Define $R_{x_1,x_2}(h)$ as the rank of a value chosen
using threshold $h$ at step $3$ given the history $(x_1, x_2)$. Its expected value is
\begin{equation}
  \label{eq:def-G}
  G_{x_1,x_2}(h) := E(R_{x_1,x_2}(h)),
\end{equation}
which can be computed directly to get 
\begin{equation}
  \label{eq:formula-G}
  G_{x_1,x_2}(h) = \frac32 + h^2 - h + (2-x_1-x_2)(1-h) + \sum_{i=1}^2
  (h-x_i)_+
\end{equation}
where $y_+ = \max(y,0)$, for all $y\in\mathbb{R}$. Then the optimal
threshold $h_3(x_1,x_2)$ must be given by 
\begin{equation}
  \label{eq:h3-argminG}
  h_3(x_1,x_2) = \argmin_{h\in[0,1]} G_{x_1,x_2}(h).
\end{equation}
For each history $(x_1, x_2)$, the graph of $G_{x_1, x_2}(\cdot)$ is
composed of the reunion of three parabolae, as illustrated in
Figure~\ref{fig:graphg}.  In this Figure we read also that the
behaviour of the minimum (mainly on which of the the three parabolae it
is to be found) depends on the region of the square $[0,1]^2$ the pair
$(x_1, x_2)$ lies in, as illustrated in
Figure~\ref{fig:h3-regions}. We do not go into detail.

Similarly as in the previous section for $RP(3)$ we need to
distinguish 5 cases, and  obtain
\begin{equation}
  \label{eq:h3-explicit}
  h_3(x_1,x_2) =
  \begin{cases}
    x_{(1)} & \text{for }(x_1,x_2)\in A_1 \\
    x_{(2)} & \text{for }(x_1,x_2)\in A_2\\
    \tilde x_1 = \frac{3-(x_1+x_2)}{2} &\text{for } (x_1,x_2)\in B_1\\
    \tilde x_2 = \frac{2-(x_1+x_2)}{2} &\text{for } (x_1,x_2)\in B_2\\
    \tilde x_3 = \frac{1-(x_1+x_2)}{2} &\text{for } (x_1,x_2)\in B_3
  \end{cases}.
\end{equation}
where the $A_i$'s and $B_i$'s are shown on
Figure~\ref{fig:h3-regions}, and where $x_{(1)}$ and $x_{(2)}$ are
respectively $\min(x_1,x_2)$ and $\max(x_1,x_2)$.

\medskip
\noindent\textbf{Step 2.} Suppose that $X_1=x_1$. The optimal
threshold $h_2(x_1)$ must be such that, if $X_2=h_2(x_1)$, then the
same payoff is obtained by selecting $X_2$ or rejecting it and acting
optimally thereafter. In other words, $h_2(x_1)$ is the indifference
value for $X_2$. Consequently the threshold  $h_2(x_1)$ must be  solution to 
\begin{equation}
  \label{eq:indif-eq-2-general}
  1 + 2h_2 + \indic(h_2>x_1) = g(x_1,x_2),
\end{equation}
with $g(x_1,x_2) := G_{x_1,x_2}(h_3(x_1,x_2)).$
The decomposition of $h_3$ given in~(\ref{eq:h3-explicit}) allows us
to obtain the explicit expression of $g(x_1,x_2)$, on each of the
regions $A_1,A_2,B_1,B_2,$ and $B_3$. After some work one notices that
the optimal threshold $h_2(x_1)$ can be obtained explicitly by
discussing separately over 6 different intervals for $x_1$.

When the history is $X_1=0$, we are faced with a
RP($3$) on $\{X_2,X_3,X_4\}$. Therefore the value of $h_2(0)$ is equal
to the value of $h_1$ in a $RP(3)$,  and (see Section~\ref{sec:expl-solut-cases})
\begin{equation}
  \label{eq:h2-0}
  h_2(0) = \frac{5-\sqrt{13}}{4} =: a.
\end{equation}
Similarly, if $X_1=1$, then we find again a $RP(3)$, hence
\begin{equation}
  \label{eq:h2-1}
  h_2(1) = a.
\end{equation}
The endcases are therefore covered. 

We now study $h_2(x_1)$ for small values of $x_1$. We know that
$h_2(x_1)$ is a continuous functions of $x_1$ (see Bruss and Ferguson
1993). The
graph of $h_2$ starts at $(0,a)$ which lies in $A_2$ (because $a>1/3$)
and ends at $(1,a)$ which lies in $A_1$ (for the same reason). We
can therefore determine $h_2$ on the interval $[0,\beta_1]$ where
$\beta_1$ is the first coordinate of the intersection of the graph of
$h_2$ with one of the boundaries of the regions $B_2$ or $B_3$. For
this reason we use the expression $G_{x_1,x_2}(x_{(2)})$
in~\eqref{eq:indif-eq-2-general} and the fact that $h_2>x_1$ when we
are close to $x_1=0$. Note that it is possible that the graph of $h_2$
intersects the line $x_2=x_1$ before it reaches the border of $B_2$ or
$B_3$. We find that the graph of $h_2$ intersects first the border
between $A_2$ and $B_3$ at the point with $x$-coordinate equal to
$\beta_1 = \frac32 \sqrt{2} - 2$. Therefore,
\begin{equation}
  \label{eq:h31}
  h_2(x_1) = \frac14\left(5 - x_1 - \sqrt{x_1^2 + 6x_1 + 13}\right) =:
  h_{21}(x_1),
\end{equation}
on $[0, \beta_1]$. 

Next, on some interval $[\beta_1, \beta_2]$ with $\beta_2$ to be
determined, we consider~\eqref{eq:indif-eq-2-general} with
$g(x_1,x_2) = G_{x_1,x_2}(\tilde x_3)$ because the graph entered the
region $B_3$. The value of $\beta_2$ is either the $x$-coordinate of
the point at which the graph of $h_2$ enters a new region, or the
point at which the solution $h_2$ of~\eqref{eq:indif-eq-2-general}
stops being strictly larger than $x_1$. Therefore, on
$[\beta_1,\beta_2]$, we have
\begin{equation}
  \label{eq:h32}
  h_2(x_1) = \sqrt{8x_1 + 54} - x_1 - 7 =: h_{22}(x_1),
\end{equation}
and we can also check that $h_{21}(\beta_1) = h_{22}(\beta_1)$. We
find that the graph of $h_2$ crosses the line $x_2=x_1$ before it
reaches another region. Therefore $\beta_2$ is the solution of
$h_{22}(x_1) = x_1$, thus $\beta_2 = \frac{\sqrt{30} - 5}{2}$.

By symmetry, these arguments also apply for large values of $x_1$
(i.e.\ close to 1).  One finds easily that
\begin{equation}
  \label{eq:h_2-654}
  h_2(x_1) =
  \begin{cases}
    \dfrac32 - \dfrac14(x_1 + \sqrt{x_1^2 - 4 x_1 + 16})    &
    \text{for } x_1\in[\beta_5,1]\\
    \sqrt{12x_1 + 42} - 6 - x_1 & \text{for } x_1\in[\beta_4,\beta_5]\\
    -\dfrac{(4x_1^2 - 6x_1 + 5)}{2(x_1 - 4)} & \text{for
    } x_1 \in [\beta_3,\beta_4]
  \end{cases}
\end{equation}
where
\begin{equation}
  \label{eq:betas}
  \beta_3 = \frac{7-\sqrt{19}}{6}, \quad \beta_4 = \frac12(11 -
  3\sqrt{11}),\quad \beta_5 = \frac12 (7 - 3\sqrt{3}).
\end{equation}
The left-hand-side of~\eqref{eq:indif-eq-2-general} was equal to
$1+2h_2$ as we started at $x_1=1$ and moved to the left. At $\beta_3$,
we have $h_2(x_1)=x_1$. At this point, $h_2(x_1)$ is not strictly
lower than   $x_1$ anymore. 

Finally we need to obtain $h_2$ for intermediate values of
$x_1 \in [\beta_2,\beta_3]$; to this end we need to consider
separately the cases $x_1\in [\beta_2,1/4)$ and
$x_1\in [1/4,\beta_3]$. We get the dichotomy (i) $h_2<x_1$ then the
\textsc{lhs} of~\eqref{eq:indif-eq-2-general} is strictly smaller than
its \textsc{rhs}, (ii) $h_2>x_1$ then the \textsc{lhs}
of~\eqref{eq:indif-eq-2-general} is strictly larger than its
\textsc{rhs}.  This can be interpreted in a probabilistic way: if
$h_2$ is taken smaller than $x_1$, the expected payoff is better if we
could stop on this value (\textsc{lhs$<$rhs}), while it is a bad
choice to stop on $X_2=h_2$ if $h_2 > x_1$ since the expected payoff
is then worse than what expected if one continues the game
(\textsc{lhs$>$rhs}). From these two observations, we  conclude
that $h_2 = x_1$.

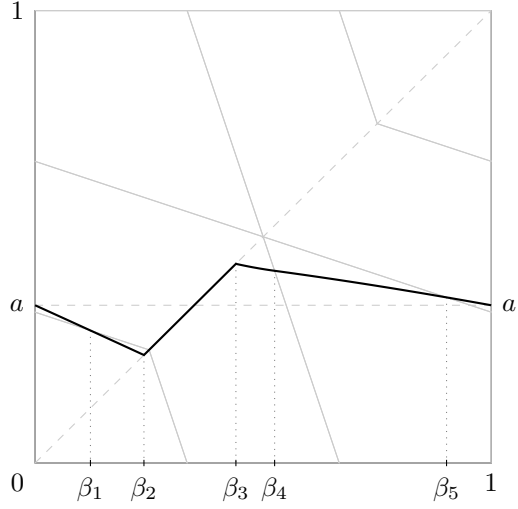
\begin{figure}[h]
  \centering
  \begin{tikzpicture}[scale=6]
    \draw[] (0,0) rectangle (1,1);
    \coordinate (q1) at (1/4,1/4);
    \coordinate (q2) at (1/2,1/2);
    \coordinate (q3) at (3/4,3/4);
    \draw[gray!40,dashed, thin] (0,0)--(1,1);
    \draw[gray!40,dashed, thin] (0,0.348612181135) --  (1,0.348612181135);
    \draw[gray!40,thin] (0,1/3) -- (q1) -- (1/3,0) -- (0,0) --
    cycle;
    \draw[gray!40,thin] (0,1/3) -- (q1) -- (1/3,0) -- (2/3,0) --
    (q2) -- (0,2/3) -- cycle;
    \draw[gray!40,thin] (0,2/3) -- (q2) -- (1/3,1) -- (0,1) --
    cycle;
    \draw[gray!40,thin] (2/3,0) -- (q2) -- (1,1/3) -- (1,0) --
    cycle;
    \draw[gray!40,thin] (1/3,1) -- (2/3,1) -- (q3) -- (1,2/3) --
    (1,1/3) -- (q2) -- cycle;
    \draw[gray!40,thin] (2/3,1)--(1,1)--(1,2/3)--(q3)--cycle;
    
    \node[left] at (0,0.348) {\footnotesize$a$};
    \node[right] at (1,0.348) {\footnotesize$a$};
    \node[left] at (0,1) {\footnotesize$1$};
    \node[below] at (1,0) {\footnotesize$1$};
    \node[anchor=north east] at (0,0) {\footnotesize$0$};
    
    \foreach \x/\y/\l in {0.12132/0.29289/1,0.23861/0.23861/2,0.44018/0.44018/3,
      0.52506/0.42482/4, 0.90192/0.36603/5}{
      \draw[dotted, thin, gray] (\x,\y)--(\x,0);
      \draw (\x,0.2pt)--(\x,-0.2pt) node[below]{\footnotesize$\beta_\l$}; }

    \draw[thick] plot[domain=0:0.12132] (\x, {0.25 * (5 - \x - sqrt(\x*\x +
      6*\x +13))}) -- 
    plot[domain=0.12132:0.2386] (\x, {sqrt(8*\x+54)-\x-7}) --
    (0.2386,0.2386) -- (0.44018,0.44018) -- 
    plot[domain=0.4402:0.5251] (\x,{-(4*\x*\x-6*\x + 5)/(2*\x-8)})-- 
    plot[domain=0.5251:0.91] (\x, {sqrt(12*\x + 42) - 6 - \x}) --
    plot[domain=0.91:1] (\x, {1.5 - (\x + sqrt(\x*\x-4*\x+16))/4});
    
  \end{tikzpicture}
  \caption{Plot of $h_2(x_1)$ for $x_1 \in [0,1] = [0, \beta_1] \cup
    [\beta_1, \beta_2] \cup [\beta_2, \beta_3] \cup [\beta_3, \beta_4]
    \cup [\beta_4, \beta_5]
    \cup [\beta_5, 1]$. Although there are 6 different expressions, it
    can be checked that  $h_2(\cdot)$ is differentiable at
    $\beta_i$ for $i\in \left\{1,  4, 5 \right\}$.}
  \label{fig:ploth2}
\end{figure}

We therefore know the expression of $h_2$ for all values of $x_1$ on
$[0,1]$; this is represented in Figure~\ref{fig:ploth2}.

\medskip
\noindent\textbf{Step 1.} The much sought-after threshold $h_1$ is
solution to
\begin{equation}
  \label{eq:indiff-h1}
  1 + 3 h_1 = g(h_1),
\end{equation}
where $g(x_1)$ is the expected rank of the selected variable if one
starts the game at step~2 with the history $X_1=x_1$ and acts
optimally thereafter.

Let us try to find a solution $h_1\in[0,\epsilon]$. The
right-hand-side of~(\ref{eq:indiff-h1}) is an integral where the
integrating variable represents the value of $X_2$; when $X_2=u\le
h_2(h_1)$, one must accept $X_2$, while one must reject $X_2=u$ if
$u>h_2(h_1)$. The behaviour when one moves on to step~3 depends on the
region the history $(h_1,u)$ lies in: $A_2$, $B_2$, or $B_3$. The
expression of $G_{h_1,u}$ will depend on this.

For the sake of concision, we will only write out the complete
expression of the integral for the smaller values of $h_1$. We thus
have 
\begin{align*}
  g(h_1) &= \int_0^{h_1} (1+2u)\,du +
  \int_{h_1}^{h_2(h_1)} (2+2u)\,du \nonumber\\
  &\phantom{=} + \int_{h_2(h_1)}^{(2-h_1)/3} G_{h_1,u}(u)\,du\\
  &\phantom{=} + \int_{(2-h_1)/3}^{1} G_{h_1,u}((1-(h_1+u))/2)\,du.
\end{align*}

The function $h_2(\cdot)$ is defined on $6$ different intervals. Thus
the need to write at least $6$ integrals in order to keep explicit
expressions around. Also look at the change in the path made
vertically through the regions $A_1,A_2,B_1,B_2,B_3$. When the regions
or the order of the regions in which we cross them changes, we must
write a separate integral. Summing things up, we need $11$ divisions
of $[0,1]$ on which the expression of the integral is each time
different. The solution to~(\ref{eq:indiff-h1}) is found on
$[\beta_2,\beta_3]$, with $\beta_2$ and $\beta_3$ defined above. The
software \textsf{Mathematica} came in handy 
for this task, yielding
\begin{align*}
  \label{eq:value-h1}
  h_1 &= \left( \tfrac{6}{1849}\sqrt{123199} -
    \tfrac{87150}{79507}\right)^{1/3} - \tfrac{846}{1849}
  \left(\tfrac{6}{1849} \sqrt{123199} -
    \tfrac{87150}{79507}\right)^{-1/3} + \tfrac{53}{43} \\
  &= 0.27502\cdots.
\end{align*}
Wrapping up we finally obtain (computations not included)
\begin{align*}
  V(4) &= -\frac{5553791}{8640} + \frac{767}{80\sqrt3} +
  \frac{2609\sqrt{11}}{216} + \frac{3281\sqrt{19}}{216} -
  \frac{59(53-\alpha_1 + \alpha_2)}{1548}\\
  &\phantom{=\:\:}+ \frac{85(53-\alpha_1+\alpha_2)^2}{44376} -
  \frac{53(53-\alpha_1+\alpha_2)^3}{2862252}
  + \frac{(53-\alpha_1+\alpha_2)^4}{11449008}\\
  &\phantom{=\:\:}+ \frac{1}{192} (842-532\sqrt3 + 31\sqrt{13} + 216
  \arcccsh(2\sqrt3)\\
  &\phantom{=\:\:}-216\arcsinh(\tfrac{3-\sqrt3}{4}) -
  \frac{2025\log(12)}{8} + \frac{2025\log(252)}{8} \\
  &\phantom{=\:\:}+\frac{1}{288}(2586985 - 779844\sqrt{11} \\
  &\phantom{=\:\:} +72900\log(\tfrac37(-1+\sqrt{11}))) - \frac{2025}{8}
  \log(17+\sqrt{19})\\
  &= 1.49329\cdots,
\end{align*}
with
\begin{align*}
  \alpha_1 &= \left(\frac{5076}{14525+43\sqrt{123199}}\right)^{1/3},\\
  \alpha_2 &= \left(6(-14525+43\sqrt{123199})\right)^{1/3}.
\end{align*}

All \textsf{Mathematica}  computations are available on Yvik Swan's
webpage.\footnote{https://sites.google.com/site/yvikswan/}
Summarizing, we have obtained the following optimal thresholds : 
\begin{gather*}
  h_1 = \left( \tfrac{6}{1849}\sqrt{123199} -
    \tfrac{87150}{79507}\right)^{1/3} - \tfrac{846}{1849}
  \left(\tfrac{6}{1849} \sqrt{123199} -
    \tfrac{87150}{79507}\right)^{-1/3} + \tfrac{53}{43},\\[4ex]
  h_2(x_1) =
  \begin{cases}
    \frac14\left(5 - x_1 - \sqrt{x_1^2 + 6x_1 + 13}\right) & \text{if
      \:} x_1 \in[0, \frac32 \sqrt{2} - 2]\\
    \sqrt{8x_1 + 54} - x_1 - 7 & \text{if \:} x_1\in [\frac32
    \sqrt{2} - 2, \frac{\sqrt{30} - 5}{2}]\\
    x_1 & \text{if \:} x_1\in[\frac{\sqrt{30} - 5}{2},
    \frac{7-\sqrt{19}}{6}]\\
    -\dfrac{(4x_1^2 - 6x_1 + 5)}{2(x_1 - 4)} & \text{if \:}
    x_1\in[\frac{7-\sqrt{19}}{6}, \frac12(11 - 3\sqrt{11})]\\
    \sqrt{12x_1 + 42} - 6 - x_1 & \text{if \:} x_1 \in [\frac12(11 -
    3\sqrt{11}), \frac12 (7 - 3\sqrt{3})]\\
    \tfrac32 - \tfrac14(x_1 + \sqrt{x_1^2 - 4 x_1 + 16}) & \text{if
      \:} x_1\in[\frac12 (7 - 3\sqrt{3}), 1]
  \end{cases},\\[4ex]
  h_3(x_1,x_2) =
  \begin{cases}
    x_{(1)} & \text{if }(x_1,x_2)\in A_1 \\
    x_{(2)} & \text{if }(x_1,x_2)\in A_2\\
    \tilde x_1 = \frac{3-(x_1+x_2)}{2} &\text{if } (x_1,x_2)\in B_1\\
    \tilde x_2 = \frac{2-(x_1+x_2)}{2} &\text{if } (x_1,x_2)\in B_2\\
    \tilde x_3 = \frac{1-(x_1+x_2)}{2} &\text{if } (x_1,x_2)\in B_3
  \end{cases},\\[4ex]
\end{gather*}
and, of course, $h_4 = 1$. 
Approximate values of the $\beta_i$'s, rounded to the 5th decimal:
\begin{gather*}
  \beta_1=0.12132,   \beta_2=0.23861,   \beta_3=0.44018,
  \beta_4=0.52506, 
  \beta_5=0.90192.
\end{gather*}

\section*{Acknowledgments}

This note was written after a conference held in Brussels on September
10 and 11, 2015 in honor of  Prof. F. T. Bruss who was supervisor of
both the authors' PHDs. His enthusiasm for his craft is an inspiration
for both of us and we thank him for his guidance during those
important early  moments of our career.

\end{document}